\pdfoutput=1 
\documentclass[12pt]{article}
\usepackage{amsmath,amsfonts,amssymb,url,color,epsfig,graphics}
\usepackage[latin1]{inputenc}

\newcommand{\R}{{\mathbb{R}}}

\newcommand{\C}{{\mathbb{C}}}
\newcommand{\Z}{{\mathbb{Z}}}
\newcommand{\N}{{\mathbb{N}}}

\def\ha{\frac{1}{2}}

\def\pa{\partial}
\def\ra{\rightarrow}
\def\ga{\alpha}

\def\gd{\delta}

\def\gg{\gamma}

\def\gl{\lambda}

\def\gs{\sigma}

\def\OPD{pseudo-differential operator}
\def\nghbd{neighbourhood~}

\newtheorem{theo}{Theorem}[section]

\newtheorem{conj}{Conjecture}[section]

\begin{document}

\title{Classical and Quantum mechanics on  3D contact manifolds\\
{\normalsize dedicated to  Victor, teacher and friend.}}

\author{Yves  Colin de Verdi\`ere\footnote{Universit\'e Grenoble-Alpes,
Institut Fourier,
 Unit{\'e} mixte
 de recherche CNRS-UGA 5582,
 BP 74, 38402-Saint Martin d'H\`eres Cedex (France);
{\color{blue} {\tt yves.colin-de-verdiere@univ-grenoble-alpes.fr}}}}

\maketitle

\abstract
In this survey paper, I  describe some aspects of the dynamics and the spectral theory of sub-Riemannian
3D contact manifolds. We use Toeplitz quantization of the characteristic cone
as introduced by Louis Boutet de Monvel and Victor Guillemin. We also discuss trace formulae following our work as well as the Duistermaat-Guillemin trace formula. 
\section{Introduction}
The goal of our work with Luc Hillairet  and Emmanuel Trélat,  in particular in \cite{C-H-T-18},
  was to see if we can extend what is known for spectral asymptotics of the
Laplace operator on a Riemannian manifold to sub-Riemannian (``sR'' in  what follows) manifolds, in particular concerning
\begin{itemize}
\item 
Trace formulae relating
the spectrum of  the Laplace operator to the lengths of periodic geodesics (see \cite{CdV-73,D-G-75} and the survey \cite{CdV-07}).
\item Quantum limits and quantum ergodicity (Schnirelman theorem, see \cite{CdV-85} and the excellent review \cite{Dy-21}).
\item Approximation of eigenfunctions by the construction of  quasi-modes, i.e. approximate solutions
of the eigenvalue equation,  supported by invariant sets of the
geodesic flow (see \cite{CdV-77,Ze-17}).
\end{itemize}
What I find  nice with sR geometry is that we have to take into account a topological set of data, namely the distribution
(a sub-bundle of the tangent bundle),  and also the metric data which allow to define the distance, the 
geodesics and the Laplace operator.
In this review, we will only speak of the case of 3D contact manifold going a little   beyond the paper \cite{C-H-T-18}. 
While starting our project, we discovered new things about the sR geodesic flow, namely  the important role played by the
Reeb vector field as a way to ``compactify'' the geodesic flow on the unit cotangent bundle whose fibers are not  compact (in contrast with the Riemannian case). We tried also to study interesting
examples like magnetic Laplacians  and the so-called Liouville Laplacians  where the sR manifold $M$ is the unit cotangent bundle 
of a Riemannian surface  and the distribution is the kernel of the Liouville form restricted to $M$.
 Special cases  of these examples linked 
to surfaces of constant curvature turn out to be ``integrable'' in a weak sense.  

\section{The setup} 

Let us consider  a 3D  smooth closed manifold $M$ equipped with a smooth contact distribution $D\subset TM $ (assumed to be oriented), a smooth metric $g$ on $D$ and
a smooth density $|dx|$.
Recall that $D$ is contact if there exists a  non vanishing real 1-form $\ga $ with $D=\ker \ga $ so that  the 3-form
$\ga \wedge d\ga $ is a volume form. 
Such a set of data $(M,D,g,|dx|)$ defines  an sR manifold with a volume form.
To such an sR  manifold, we associate the following objects:
\begin{enumerate}
\item The cometric $g^\star :T^\star M \ra \R^+ $  is defined in local coordinates by
\[ g^\star (x,\xi)=\| \xi_{|D_x} \|^2_{g(x)} \]
\item The geodesic flow which is the flow of the Hamiltonian vector field $X_g$ of $\ha g^\star $. When restricted
to the unit cotangent bundle $U^\star M:=\{ (x,\xi)\in T^\star M\ |\ g^\star (x,\xi)=1 \} $
the integral curves project onto $M$ as geodesics with unit speed and, conversely, any geodesic  is the projection of such an integral curve. 
\item The Laplacian $\Delta _{sR} $
  which is the Friedrichs extension on $L^2 (M,|dx|)$ of the quadratic form
  $q(f)=\int_M g^\star (df) |dx| $.  The self-adjoint second order differential operator $\Delta _{sR} $
  can be written locally as $\Delta =X^\star X + Y^\star Y =-X^2 -Y^2 +{\rm l.o.t} $
  where $(X,Y)$ is a local smooth  orthonormal  frame of $D$ and $X^\star $ and $Y^\star $ are the adjoints of $X$ and $Y$ with respect to $|dx|$. 
The operator $\Delta_{sR} $ is sub-elliptic (a well known result due to Lars H\"ormander) and hence has a discrete spectrum $(\gl_j),~j\in \N$, with an o.n.b. $(\phi_j),~j\in \N, $
of eigenfunctions $(\phi_j),~j\in \N$ of $L^2$. 
The principal symbol of this Laplacian is the co-metric. 
\item The canonical contact 1-form $\ga _g $ which is defined by $\ker \ga _g =D$ 
and $(d\ga_g )_{|D}= v_g $ where $v_g$ is the volume form on $D$ induced by the metric and the orientation. 
\item The characteristic manifold $\Sigma :=D^\perp =(g^\star)^{-1}(0) $ which is a 4D symplectic subcone  of $T^\star M\setminus 0 $. 
\item The  Reeb vector field $R$ of the form $\ga_g$ (ie $\ga_g (R)=1,~\iota (R)d\ga_g=0$) which is the projection onto
$M$ of the Hamiltonian vector field of 
$\rho:\Sigma \ra \R $ defined by 
$\rho (s\ga_g )=s$.  
\end{enumerate}
The main object of this review  is to describe asymptotic properties of the geodesic flow and of the
spectral data of the  sR Laplacian. 
The geodesics of large momenta spiral around the Reeb flow. This leads to the existence of infinitely many  periodic geodesics spiraling around a
generic closed Reeb orbit.
Concerning the Laplace operator, most eigenfunctions concentrate microlocally on the characteristic manifold.
Here a natural  object  is the quantization of the Reeb Hamiltonian as a Toeplitz operator
``\`a  la Boutet de Monvel/Guillemin'' (\cite{B-G-81}). We recover band spectra which we call ``Landau bands'': they are indeed Landau levels in
some  magnetic  examples. 
The starting point of this presentation is works in collaboration with Luc Hillairet (Orl\'eans) and Emmanuel Tr\'elat (Paris).
There  are no fundamentally   new results in the present paper, but some  new definitions, examples  and conjectures.
The  conjecture that I propose (see Section \ref{sec:conj}) is the following one
\begin{conj}
 The periods of the Reeb orbits are spectral invariants of the sR Laplacian.
\end{conj}

\section{Example 0: Heisenberg quotients}
This is the most basic example (see \cite{C-H-T-18}, sec. 3.1).
We consider  the presentation of $H^3$ as $\R^3$ equipped with the group law
$(x,y,z)\star (x',y' ,z')=(x+x',y+y', z+z'-xy')$. 
 We choose the subgroup $\Gamma :=\{ (x,y,z)\ | \ (x,y)\in (\sqrt{2 \pi }\Z )^2,~ z\in 2\pi \Z \}$.
Our sR manifold is then $\R^3/\Gamma $ with 
the orthonormal basis  for $D$ given by 
\[ X=\pa_x,\  Y=\pa_y - x \pa _z \]
The spectrum of $\Delta =-(X^2+Y^2)$ is then explicitly computable:
one gets the union of the  eigenvalues of the flat torus $\R^2/\sqrt{2\pi }\Z^2$ 
and the set of integers $m (2l+1),~ m=1, \cdots,~l=0,\cdots $ with multiplicities $2m$. 
Note that the multiplicities are very high.

The lengths spectrum (the set of lengths of closed geodesics) is the set of $2\pi \sqrt{2n},~n\in \N$. 

Note that $\ga_g=dz+xdy $ and the Reeb vector field is $\pa_z$ which is a Killing vector field. 

\section{Example 1: magnetic fields  over a Riemannian  surface}\label{ex:B}
Let  $\pi: M \ra X $ be a principal $S^1_\theta $-bundle on an oriented Riemannian  surface $(X,h)$.
We assume that this bundle is equipped with an Hermitian  connection $\nabla $ whose horizontal distribution is our $D$. If the curvature of
the connection does not vanish, the distribution $D$ is contact. We take for $g$ the pull-back on $D$ of the metric $h$ by
$\pi$. The curvature of $\nabla $ is a 2-form $B $ (the magnetic field)  and  one introduces the
magnetic scalar $b=B/ dx_h $ where $dx_h$ is the Riemannian volume form of $X$.
The sR metric is invariant by the $S^1$ action, this gives an invariant momentum $e:T^\star M \ra \R $ which is the principal symbol of $-i\pa _\theta $. 
The geodesics of $(M,D,g)$ with momentum $e$ project onto the trajectories on $X$ with the magnetic fields $b$ and electric charge $e$.

The Reeb flow is
\[ R=  b \partial _\theta - \vec{b} \]
where  $\vec{b}$ is the horizontal lift of the Hamiltonian vector field of $b$ w.r. to the symplectic form $B$ on $X$.
We define the Laplacian using the volume form
$dx_g= |d\theta \wedge \pi^\star dx_h |$.
Then $\Delta $ commutes with the $S^1$-action and $L^2(M, dx_g)$ splits into a direct sum
$\oplus _{n\in \Z} H_n $ where $H_n $ is unitarily equivalent to the Schr\"odinger operator on $X$ with magnetic field $nB$.

\section{Example 2: Liouville form on the unit cotangent bundle of a Riemann surface}

Again $(X,h)$ is a Riemannian  surface and  $M$ is the unit cotangent bundle of $X$.
The distribution $D$ is $D=\ker \gl $ where 
$\gl $ is the restriction to $M$ of the Liouville 1-form $\xi dx +\eta dy $. 
We take on $D$  any metric so that $\ga_g =\gl$.
Then the Reeb vector field is the geodesic flow of $h$. 

In particular, the case of hyperbolic surfaces is of special interest.
This example is studied in \cite{C-H-W-?}. 
Using the representations of ${\rm SL}_2(\R)$, we reduce the computation  of the spectrum and periodic geodesics to 
1D-problems.

\section{Example 3: Jacobi metric for the sR Kepler problem}

Following \cite{Sh-21}, we take the metric 
 $g= D^{-\ha}g_0 $ on $\R^3\setminus 0$ where  $g_0$ is the Heisenberg metric
and $D=(x^2+y^2)^2 + 16 z^2$. This metric is the Jacobi metric for a sR Kepler problem at energy $0$.
The metric $g$ is invariant by the dilations
$\gd_\gl: (x,y,z) \ra (\gl x, \gl y ,\gl^2 z)$. 
It admits compact quotients by the groups
generated by $\gd _{\gl_0}$ for some  $\gl_0>1$. Hence the geodesic flow is complete and the Laplacian is essentially self-adjoint.

\section{Example 4:  boundary of complex domains}
If $\Omega $ is a smooth  domain in $\C^2$, we consider, on $M:=\pa \Omega $,  the distribution $D=TM\cap iTM $.
If $\Omega $ is strictly pseudo-convex, $(M,D)$ is a contact manifold. We can take the metric induced by the Euclidean
metric on $\C^2$.

Another example after Louis Boutet de Monvel \cite{BdM-80}: 
let $Z\subset \C^N $ be a complex subcone of complex dimension $2$
(for example defined as the zero set of complex valued  homogeneous polynomials),  smooth outside $0$,
 and $B$ the unit ball of $\C^N$.
The 3D manifold  $M = \pa (Z \cap B)$ is an  $S^1$-bundle over the projective complex curve $Z\setminus 0 / \C \setminus 0 $.  
The form $\ga =\sum  _j \Im (z_j d\bar{z}_j )$ is contact on $M$.
The Reeb flow $R$ of $\ga$ is $2\pi $ periodic.  A convenient choice of $g$ has the Reeb flow $R$. We call Such a manifold a
Zoll-Reeb sR manifold.

\section{Classical Birkhoff normal forms } \label{sec:gnf}

We will assume for simplicity that the fiber bundle $D\ra M$ is topologically trivial.
This holds for the magnetic sR if $X$ is a torus.
This holds for the Liouville sR if the surface $X$ is orientable.
This holds also if $M$ is a \nghbd  of a periodic  Reeb orbit. 
Then
 \begin{theo}
 \label{theo:BNF}
  There exists an homogeneous canonical transformation
$ \chi :C\ra C'$ with $C$  a conic \nghbd of $\Sigma $ in $T^\star M\setminus 0$
and $C'$ a conic \nghbd of $\Sigma \times 0  $ in $\Sigma \times \R^2$, with $\R^2_{u,v} $ equipped with the symplectic form
$dv\wedge du $ and the cone structure
$\gl.(u,v)=(\sqrt{\gl}u,\sqrt{\gl}v)$, 
so that $\chi _{|\Sigma }={\rm Id}\times 0$
and
\[ g^\star \circ \chi^{-1} = \sum _{j=1}^\infty \rho _j(\sigma ) I^j + O \left((I/\rho \right)^{\infty})  \]
with
$\rho_j :\Sigma \setminus 0 \ra \R $ homogeneous of degree $2-j$, $\rho_1=|\rho | $  with $\rho $ the Reeb Hamiltonian 
and $I=u^2 +v^2$.  The function  $\rho_2$ is uniquely defined modulo  Lie derivatives w.r. to Reeb.

\end{theo}

This is proved in Section 5 of \cite{C-H-T-18}.

\section{Spiraling of the sR geodesics around Reeb orbits}

The goal of this section is to explain the following fact:
given $x_0 \in M$ and $v_0 \in D$ of length $1$ for the metric $g(x_0)$, there exists a 1-parameter family of geodesics with these Cauchy data at time $0$.
They are associated to initial momenta whose component vanishing on $D$ is not fixed. When this transverse momentum tends to
$\infty $, these geodesics will spiral more and more around a Reeb orbit like helices with small radii.
See \cite{C-H-T-21} for more details. 

\subsection{A simple  Hamiltonian}
Let us assume that our Hamiltonian is 
$H_0=\ha \rho I $ on $\Sigma \times \R^2 $ with $\rho $ the Reeb Hamiltonian. 
Then the Hamiltonian vector field is
\[ \vec{H_0}= \ha I \vec{\rho}+ \rho \pa _\theta \]
The Poisson bracket $\{ \rho, I\} $ vanishes, hence $\rho $ and $I$ are first integrals of the motion. 
The dynamics  can be integrated as follows:
\[ \Phi_t (\gs_0, u+iv )= (\phi_{It/2}(\gs _0), (u+iv)e^{i\rho (\gs_0) t}) \]
where $\Phi_t$ is the  flow of $H_0$, $\phi_t$ the flow of $\rho$ (the Reeb flow). 
If  we fix the energy $H_0=1$, we have
\[ \Phi_t (\gs_0, u+iv )= (\phi_{It/2}(\gs _0), (u+iv)e^{i t/I}) \]
As $I$ is small, $\rho $ is large and we get a spiraling flow around the Reeb orbits. Note that this Hamiltonian is exactly the
Heisenberg one. In particular, there exists closed geodesics $\gg_k,\ k\in \N ,$ spiraling around any periodic Reeb orbit of period $T_0$ of lengths
$l_k=2\sqrt{\pi k T_0}$. 
\subsection{Spiraling}
Let us choose an orthonormal frame $(X,Y)$ of $D$ in some tubular \nghbd of  a Reeb orbit $\Gamma$ defined on some interval $t\in [0,T]$. Denote by $Z=[X,Y]$.
Moreover, we can deduce from the Birkhoff normal form  the existence of a well defined parallel transport of vectors in $D$ along the Reeb flow.
\begin{theo}\label{theo:spiraling}
Let $q_0\in M$ be arbitrary, and let $(q_0,p_0)\in T^\star_{q_0} M$ be the Cauchy data of a geodesic $t\mapsto\gamma (t)$ starting at $q_0$
with unit speed $\dot{\gamma }(0)=X_0 \in D(q_0)$. We assume that $p_0 \ra \infty $ (large initial momentum) and denote by $h_0=p_0(Z)\ra \infty $.

Then, there exists a point $Q_0=Q_0(q_0,p_0)\in M$ close to $q_0$, and a vector $Y_0 \in D(Q_0)$ close to $X_0$, such that,
denoting  by $\Gamma (\tau )= \mathcal{R}_{\tau }(Q_0)$ the Reeb orbit of $Q_0$, and by $Y(t)$ the   parallel transport of $Y_0$ along $\Gamma $, 
we have, using the complex structure on $D$, for $t=\mathrm{O}(h_0)$, 
$$
\gamma (t)=  \Gamma (J_0t/2)-i J_0 e^{it/J_0}Y(J_0t/2)) + \mathrm{O}(J_0^2) $$

with $J_0=h_0^{-1}+ \mathrm{O}(h_0^{-3}) $.
\end{theo}
In words, the sR geodesic $\gamma $ spirals along a Reeb orbit with a slow speed $\sim 1/2h_0$ along that orbit and a fast angular speed
$\sim h_0$ transversally. 

\subsection{Periodic geodesics around a generic periodic Reeb orbit}

In the paper \cite{CdV-22a}, I proved the following
\begin{theo}\label{theo:per}  If $\Gamma $ is a non degenerate periodic orbit of the Reeb flow of period $T_0>0$,
there exist infinitely many periodic sR geodesics
$(\gamma_k),~k\geq k_0$, accumulating on $\Gamma$ as $k\ra +\infty$, 
whose  lengths admit a full asymptotic expansion
\[ L_{k}=2\sqrt{\pi T_0}k^\ha  + \sum _{j=1}^\infty a_{j} k^{-j/2} +O\left(k^{-\infty}\right) \]
as $k\ra +\infty $. 
\end{theo}

\subsection{Periodic geodesics in the weakly integrable case}\label{sec:w}

We say that the geodesic flow is {\it weakly integrable}  if the BNF converges, ie we can write
$g^\star =F(\sigma ,I )$ in some conic \nghbd  of $\Sigma $ with $F$ admitting an expansion in powers of $I$
as given in Theorem \ref{theo:BNF}.
Let us sketch a proof of Theorem \ref{theo:per} in this case.

 First, if  $I$ is  small enough,  there exists a closed orbit $\Gamma _I $, 
of the Hamiltonian $\rho +\rho_2 I + \cdots $ of period
$T(I)=T_0 -AI + \cdots $ contained in $H(\sigma ,I)=1$.  
One then consider  at the return map of the angles. This gives
\[ \int_0^{T(I)}\left(\rho +2I \rho_2 +\cdots \right) dt =2k\pi \]
One then use the fact that $H=1$ and concludes by eliminating $I$. 

Note that this asymptotic formula is exact in the examples  of Heisenberg ($T_0=2\pi $)
and $S^3$ ($T_0=\pi $) (see \cite{K-V-19}). 

\section{Weyl measures, QL and QE for the sR Laplacian}
\subsection{Weyl}
The Weyl formula is given by: 
\begin{theo} 
If $N(\gl )=\# \{ j \ | \ \gl_j \leq \gl \}$, we have, as $\gl \ra +\infty$, 
\[ N(\gl )\sim \frac{ \int_M |\ga_g \wedge d\ga_g |}{32} \gl^2 \]
\end{theo}
Note that the exponent $2$ of $\gl$ is larger than the exponent $3/2$ of  the Riemannian case.
The smooth measure  $|\ga_g \wedge d\ga_g |$ is called the Popp measure. 
Note that  the measure $\mu $ which is used in order to define the Laplacian does not need to be the Popp measure.  
The measure $\mu $ plays a very minor role in the spectral asymptotics. This is because, for any pair  $\mu, \ \mu_0$ of densities, 
$\Delta_{g,\mu} $ is unitary equivalent to $\Delta_{g,\mu_0} +V $ for some smooth potential $V$. 
Hence, one gets, for example, by using the minimax  principle, that
\[ \exists C>0 {\rm ~so~ that}, \forall j \in \N, ~|\gl_j (\Delta_{g,\mu})-\gl_j (\Delta_{g,\mu_0})|\leq C \] 

The volume $\int_M |\ga_g \wedge d\ga_g |$ which is a spectral invariant corresponds to  the inverse of 
Arnold's asymptotic linking number for the Reeb flow if $M=S^3$
\cite{Ar-86}.

There is a microlocal version of the Weyl law, namely
\begin{theo} 
If $A$ is self-adjoint pseudo-differential operator of degree $0$ whose principal symbol $a:T^\star M\setminus 0\ra \R $ is  homogeneous 
of degree $0$ and is identified with a function on the sphere bundle $S(T^\star M)$,  
we have
\[ \lim_{\gl \ra +\infty }\frac{1}{N(\gl) } \sum _{\gl_j \leq \gl}
\langle A\phi_j |\phi_j \rangle= \int _{S(\Sigma )} a dL \]
where $dL $ is the unique probability measure on $S(\Sigma )$ which is invariant by antipody and
whose direct image by the projection onto $M$ is the probability measure $|\ga_g \wedge d\ga_g|/\int_M  | \ga_g \wedge d\ga_g |$.
\end{theo}
Both theorems are proved in \cite{C-H-T-18}. The first one is classical, but we provided a new proof in that paper.

\subsection{QL and QE}
Let us recall what are {\it quantum limits} (in short QL's): a QL is a probability measure $dm $ on $S^\star M:=S(T^\star M)$ such
that there exists a sequence of eigenfunctions $\phi_{j_k},~k\in \N, $
 of our Laplacian such that for any self-adjoint \OPD ~$A$ of degree $0$ with 
homogeneous principal symbol $a\in C^\infty (S^\star M)$, on has
\[ \lim_{k\ra +\infty }\langle A\phi_{j_k}|\phi_{j_k} \rangle = \int _{S\star M} a dm \]

The eigenbasis $\phi_j, j\in \N, $ is said to satisfy {\it  Quantum Ergodicity} (in short QE) with a probability measure
$dE$ on $S^\star M$  if there exists a subsequence  $(j_k),k\in \N,  $ of density one w.r. to the Weyl law, 
admitting $dL $ as QL. 

We have the following  results  which shows the prominent role of the Reeb vector field in the spectral asymptotics
(\cite{C-H-T-18}, Theorems A and B):
\begin{theo}
Let us  decompose the sphere bundle $S(T^\star M \setminus 0)$ as the disjoint union of the
unit bundle $U^\star M:=\{ g^\star =1\} $ and the sphere bundle of the characteristic manifold $S\Sigma $.
\begin{enumerate}
\item Any QL $\mu $ (a probability measure on $S(T^\star M \setminus 0)$), can be uniquely
written as the sum $\mu=\mu_0+\mu _\infty $ where  $\mu_\infty $ is supported by $S\Sigma $ and is invariant under the Reeb flow, while
 $\mu_0 (S\Sigma )=0$ and $\mu_0$ is invariant under  the geodesic flow.
\item If $(\phi_j)$ is an ONB of eigenfunctions, there exists a subsequence $(\phi_{j_k})$ of density $1$,
so that all corresponding QL's are supported on  $S\Sigma $ (and hence invariant by Reeb).
\item If the Reeb flow is ergodic, then we have QE for any real eigenbasis with the limit measure the measure $dL$ on $S\Sigma $. 
\end{enumerate}
\end{theo}

\section{Toeplitz quantization of the Reeb Hamiltonian and Landau levels}

Recall that, if $\Sigma $ is a symplectic sub-cone of $T^\star M$, one can associate to it an Hilbert space ${\cal H}\subset L^2$
of functions whose wavefront set is included in $\Sigma $ and an algebra of operators which obey to the usual rules of the pseudo-differential calculus
 where the   symbols are functions on $\Sigma $ (see \cite{B-G-81}).
In particular in the case  of boundaries of complex domains, one recovers the original definition of Toeplitz operators.

We have the following normal form:
\begin{theo}
Assuming  that $D$ is a trivial bundle,
we can use a FIO associated to the canonical transform $\chi $ defined in Section \ref{sec:gnf}
to transform $\Delta $ 
into 
\[ \Delta _0 = \sum _{j=1}^\infty R_j \otimes \Omega ^j + R_\infty \]
where the 
$R_j $ are Toeplitz operators on $\Sigma $ of degree $1-j$, $R_0$ is elliptic with symbol $|\rho |$,
$\Omega $ is an harmonic oscillator on $\R^2 $ and $R_\infty $ is smoothing along $\Sigma $. 
\end{theo}
The proof is a standard extension of the classical Birkhoff normal form using Fourier integral operators. 

It follows that we have, for each value of $l\in \N$ a sequence of eigenvalues of $\Delta $ 
which are the eigenvalues
of the Toeplitz operator
\[ \Delta _l:= \sum_{j=1}^\infty (2l+1)^j R_j \] 
modulo  a fast decaying sequence.

We call the spectrum of $\Delta _l $  a ``Landau band'' because in the case of constant magnetic field on surfaces (see Section
\ref{ex:B})
it is the union of the $l$-th Landau clusters for all magnetic field $kB,~k\in \Z\setminus 0 $.

Note that the approximate eigenvalues given by the normal form are ``almost all'' eigenvalues:
we have
\[ N_l(\lambda ) \sim \frac{\gl^2}{(2\pi (2l+1))^2 }\int _{|\rho |\leq 1} dL_\Sigma \]
with $dL_\Sigma $ the Liouville measure on $\Sigma $. 
This gives
\[ \sum _l N_l(\gl) \sim \frac{\gl^2}{32} \int_M |\ga_g\wedge d\ga _g | \]
which fits with the Weyl formula. 

\section{$\Gamma \backslash PSL_2(\R)$}\label{sec:hyp}
Here $\Gamma $ is a lattice in $G=PSL_2(\R)$.
We will look at operators invariant under  left translations. Their symbols are functions on the dual of the Lie algebra.
The Lie algebra is the 3D space of trace free $2\times 2$ matrices
\[ M(x,y,z):=\left( \begin{matrix} z & x \\ y & -z \end{matrix} \right) \] 
We write $M=z A+ xX  + yY $.
The Casimir operator
is $\square = -A^2 - 2(XY + YX)$. 
The Liouville Laplacian is
\[ \Delta _L =-(X^2+Y^2 ) \]
The magnetic Laplacian is
\[ \Delta _B = -A^2 - (X+Y)^2 \]
Their principal symbols are
$c = \zeta^2 +4\xi \eta $,
$l= \xi^2 +\eta ^2$ and $b= \zeta^2 +(\xi +\eta )^2$ respectively. 
The characteristic cones  are
$\Sigma _l=\{\xi=\eta =0 \} \subset  \{ c>0 \} $
and
$\Sigma _b = \{ \zeta =0, \xi+\eta =0 \} \subset \{ c<0 \}  $.
The co-adjoint orbits lying  in $c>0$ support  the principal series of irreducible representations  (H1 hyperbolo\"ids), while the
orbits lying  in $c<0 $ 
  support  the discrete  series of irreps ($H2$  hyperbolo\"ids). 

The calculation of the action of these operators on irreducible representations is the subject of  \cite{C-H-W-?}.
The spectrum of $\Delta _B $ is described in \cite{Ch-20}.
Note that the magnetic Hamiltonian is Liouville integrable while the Liouville one is
only weakly integrable (see Section \ref{sec:w}) thanks to the Euler equations. 

\section{Traces}        
\subsection{Wave traces}\label{sec:wt}
Richard Melrose  proved in the paper \cite{Me-84}   that the Duistermaat-Guillemin trace formula applies for 
the singularities of ${\rm Trace ( exp}(-it\sqrt{\Delta}))$
outside $t=0$. I gave in \cite{CdV-22b} a simpler proof of this result. 

The wave traces of the $\Delta_l$'s: the $\Delta _l$'s  are self-adjoint elliptic Toeplitz operators of degree $1$
to which the Theorems 9 and 10
of \cite{B-G-81}  apply. The corresponding closed orbits are the Reeb orbits. 
These theorems  say in particular that the singular support of the distributions
$Z_l(t):={\rm trace}({\rm exp}(it \sqrt{\Delta _l }))$ is contained in the set of periods of the Reeb flow
divided by $(2l+1)$. In fact, under some
genericity assumption on the Reeb flow the two sets are the same.
Summing with respect to $l$ gives a dense set of singularities. It would be nice to say more on the
precise structure of these singularities. 

\subsection{Schr\"odinger trace}\label{sec:tr-h3}
The Heisenberg case:  
for $\Re (z) >0$, one defines
$Z(z):=\sum_{j=1}^\infty  e^{-z\gl_j } $
For the flat Heisenberg, one gets
\begin{multline*}  Z_o (z)= \sum _{m=1}^\infty2m \sum _{n=0}^\infty  e^{-(2n+1)m z }=
\ha \sum _{m=1}^\infty \frac{m}{\sinh m z } \\
 =\frac{1}{2z} \sum _{m\in \Z} \frac{mz}{\sinh m z }-1 \end{multline*} 
We will   apply the Poisson summation formula.  
The Fourier transform of $z/\sinh z $ is
$\frac{\pi^2}{1+\cosh \pi z }$.
We get
\[  Z_o (z)= \frac{\pi^2}{4z^2}-\frac{1}{2z}+ \frac{\pi^2}{z^2 }\sum_{l=1}^\infty \frac{1}{1 + \cosh (2\pi^2 l/z )} \]
One recovers the Weyl law, $Z_o(t)\sim \pi^2/8t^2$.

As $z\ra 0^+$, one gets exponential corrections as in \cite{CdV-73}
\[ {\pi^2}\sum_{l=1}^\infty e^{-2\pi^2 l/z}\]
Each exponent identifies with $L^2/4$, hence
 one gets also the lengths of periodic geodesics, namely the set
of $2\pi \sqrt{2l},l\in \N $.

Let  us look at the boundary values as $\Re (z) \ra 0^+$.
There are infinitely many ``poles'', namely the zeroes of
$1+\cosh 2\pi^2 l/i\tau $. The poles are
$z_l= \frac{2\pi il}{2k+1}$. This corresponds to the periodic orbits of $(2k+1)$-times Reeb as expected.
This is a dense set in the boundary, there exist no meromorphic extensions.

\section{A conjecture}\label{sec:conj}

From what we know, I propose the following   conjecture:

{\bf  ``the periods of the Reeb flow  are generically spectral invariants of the sR Laplacian''. }

There are two heuristic arguments for that:
\begin{enumerate}
\item Using  Theorem \ref{theo:per},  one    recover the lengths of the periodic geodesics. Then the 
 asymptotics of the lengths of closed geodesics accumulating around a Reeb periodic orbit involve  the Reeb periods (see Theorem \ref{theo:per}).  
\item 
The second argument follows from the Boutet de Monvel-Guillemin trace formula appplied to each of the $\Delta_l$'s
(see Section \ref{sec:wt}) whose righthandsides involve the Reeb periods. 

Of course, it would be even nicer to extend the Schrödinger  trace formula from Section \ref{sec:tr-h3} to a more general case. 
\end{enumerate}

\bibliographystyle{plain}

\begin{thebibliography}{}

\end{thebibliography}


\begin{thebibliography}{99}

\bibitem[Ar-86]{Ar-86}
Vladimir Arnol'd,
\textit{The asymptotic Hopf invariant and its applications},
Selecta Math. Soviet. {\bf 5}:327--345 (1986). 


\bibitem[BdM-80]{BdM-80}
Louis Boutet de Monvel. 
{\it Nombre de valeurs propres d'un  op\'erateur elliptique
et polyn\^ome de Hilbert-Samuel.}
S\'eminaire N. Bourbaki, exp. no 532, p. 120-131 (1980).

\bibitem[B-G-81]{B-G-81}
Louis Boutet de Monvel \& Victor Guillemin.
{\it The spectral Theory of Toeplitz operators.}
Ann. Math. Studies, Princeton (1981).

\bibitem[Ch-20]{Ch-20} Laurent Charles.
{\it Landau levels on a compact manifold.}
ArXiv:2012.14190 (2020). 

\bibitem[CdV-73]{CdV-73}
 Yves Colin de Verdière.
{\it  Spectre du Laplacien et longueurs des géodésiques périodiques II.}
 Compositio Mathematica {\bf 27 (2)}:159--184 (1973). 

\bibitem[CdV-77]{CdV-77}
Yves Colin de Verdière.
{\it  Quasi-modes sur les variétés Riemanniennes.}
 Inventiones {\bf 43}:15--52 (1977).

\bibitem[CdV-85]{CdV-85}
Yves Colin de Verdière.
{\it  Ergodicité et fonctions propres du Laplacien.}
 Commun. Math. Phys. {\bf 102}:497--502 (1985). 


\bibitem[CdV-07]{CdV-07}
Yves Colin de Verdière.
{\it  Spectrum of the Laplace operator and periodic geodesics: thirty years after.}
 Ann. Institut Fourier {\bf 57}:2429--2463 (2007).

\bibitem[CdV-22a]{CdV-22a}Yves Colin de Verdière.
{\it Periodic geodesics for contact sub-Riemannian 3D manifolds.}
ArXiv:2202.13743 (2022).



\bibitem[CdV-22b]{CdV-22b}Yves Colin de Verdière.
{\it A proof of a Melrose's trace formula.}
ArXiv:2205.08744 (2022).

\bibitem[C-H-T-18]{C-H-T-18}
Yves Colin de Verdi\`ere, Luc Hillairet \&   Emmanuel Tr\'elat.
{\it Spectral asymptotics for sub-Riemannian Laplacians I: Quantum ergodicity and quantum limits in the 3D contact case},
Duke Math.
J., {\bf 167(1)}:109--174 (2018).

\bibitem[C-H-T-21]{C-H-T-21}
Yves Colin de Verdi\`ere, Luc Hillairet \&   Emmanuel Tr\'elat.
{\it Spiraling of sub-Riemannian geodesics around the Reeb flow in the 3D contact case.} 
ArXiv:2102.12741 (2021). 


\bibitem[C-H-W-?]{C-H-W-?}
Yves  Colin de Verdi\`ere, Joachim Hilgert \& Tobias Weich. 
{\it Irreducible representations of $SL_2(\R)$ and the Peyresq's operators. }
In preparation.

\bibitem[D-G-75]{D-G-75} Hans Duistermaat \& Victor Guillemin.
  The spectrum of positive elliptic operators and periodic bicharacteristics.
  {\it Invent. Math.} {\bf 29}:39--79 (1975). 

\bibitem[Dy-21]{Dy-21} Semyon Dyatlov.
{\it Around quantum ergodicity.}
ArXiv:2103.08093  (2021). 



\bibitem[K-V-19]{K-V-19}
David Klapheck \& Michael  VanValkenburgh. {\it The length spectrum
of the sub-Riemannian three-sphere}.  Involve {\bf 12}:45--61 (2019).



\bibitem[Me-84]{Me-84}
Richard B. Melrose.
\textit{The wave equation for a hypoelliptic operator with symplectic characteristics of codimension two},
J. Analyse Math. {\bf 44}:134--182 (1984).



\bibitem[Mo-02]{Mo-02}
Richard  Montgomery.
{\it A tour of subriemannian geometries, their geodesics and applications},
Mathematical Surveys and Monographs {\bf 91}, American Mathematical Society, Providence, RI (2002).



\bibitem[Mo-91]{Mo-91}
Richard  Montgomery. {\it How much does a rigid body rotates? A Berry phase from the 18th century.}
Amer. J. Phys. {\bf 59 (5)}:394--398 (1991).  

\bibitem[Sh-21]{Sh-21}
Corey Shanbrom.
{\it An introduction to the Kepler-Heisenberg problem.}
ArXiv:2101.03639  (2021). 




\bibitem[Ta-86]{Ta-86}
Michael Taylor. 
{\it Non commutative harmonic analysis.}
Mathematical Surveys and Monographs {\bf 22}, American Mathematical Society, Providence, RI (1986).



\bibitem[Ze-17]{Ze-17} 
Steve Zelditch. 
{\it Eigenfunctions of the Laplacian on a Riemannian manifold.}
CBMS Regional Conference Series in Mathematics {\bf 125}
 Providence, RI: American Mathematical Society (2017). 



\end{thebibliography}

\end{document}